\documentclass[preprint,sort&compress]{elsarticle}

\usepackage{framed,multirow}

\usepackage{amssymb}
\usepackage{latexsym}

\usepackage{url}
\usepackage{xcolor}
\definecolor{newcolor}{rgb}{.8,.349,.1}

\usepackage{fullpage}

\usepackage{amscd}
\usepackage{amsmath, amsthm, mathtools}
\usepackage{multirow}
\usepackage{rotating}
\usepackage{epsfig}
\usepackage{overpic}
\usepackage{hyperref}
\theoremstyle{definition}
\usepackage{csquotes}

\newtheorem{theorem}{Theorem}[section]

\newtheorem{remark}[theorem]{Remark}

\newcommand{\ve}{\mathbf{b}}
\newcommand{\vd}{\mathbf{a}}
\newcommand{\vn}{\mathbf{n}}
\newcommand{\ep}{\varepsilon}
\newcommand{\ds}{\displaystyle}
\newcommand{\vx}{\mathbf{x}}
\newcommand{\vxi}{\boldsymbol{\xi}}
\newcommand{\vy}{\mathbf{y}}
\newcommand{\vc}{\mathbf{c}}

\newcommand{\vu}{\mathbf{u}}
\newcommand{\vw}{\mathbf{w}}
\newcommand{\vs}{\mathbf{s}}

\newcommand{\rots}{\mathbf{L}^{\ast}}
\newcommand{\rotc}{\mathbf{L}}
\newcommand{\Pdiv}{\widetilde{P}_{\rm div}}
\newcommand{\Psidiv}{\Phi_{\rm div}}
\newcommand{\Psidivt}{\widetilde{\Phi}_{\rm div}}
\newcommand{\mutrunc}{\mu_{\text{trunc}}}
\newcommand{\Qx}{Q(\vx)}
\newcommand{\Qy}{Q(\vy)}
\newcommand{\Qyj}{Q(\vy_j)}

\begin{document}

\begin{frontmatter}

\title{A stable algorithm for divergence-free radial basis functions in the flat limit}%

\author[1]{Kathryn P. Drake}
\ead{KathrynDrake@u.boisestate.edu}
\author[1]{Grady B. Wright\corref{cor1}}
\cortext[cor1]{Corresponding author}
\ead{gradywright@boisestate.edu}

\address[1]{Boise State University, Department of Mathematics, 1910 University Drive, Boise, Idaho 83725-1555, USA}


\begin{abstract}
The direct method used for calculating smooth radial basis function (RBF) interpolants in the flat limit becomes numerically unstable. The RBF-QR algorithm bypasses this ill-conditioning using a clever change of basis technique. We extend this method for computing interpolants involving matrix-valued kernels, specifically surface divergence-free RBFs on the sphere, in the flat limit. Results illustrating the effectiveness of this algorithm are presented for a divergence-free vector field on the sphere from samples at scattered points.
\end{abstract}



\end{frontmatter}
 
%
%
\section{Introduction}
{Tangential vector fields to a sphere are important in many areas of geophysical sciences, from the surface of the ocean to the ionosphere~\cite{FPLM18}.}  Often the values of these vector fields may only be known at ``scattered'' locations, e.g., from measurement taken from rawinsondes, airplanes, buoys, remote sensing devices, or from output from certain numerical models (cf.~\cite[\S 4]{Williamson:2007DynamicalCores}), and values of the field must be approximated at other locations, e.g., on a grid or mesh.  Additionally, these tangential fields may satisfy certain physical constraints, such as being surface divergence-free (div-free) or curl-free, that must be preserved in the approximation.  A radial basis function (RBF) technique was developed exactly for these applications in the papers~\cite{NarcWardWright,FuselierWright2009}.  The idea is to construct a positive definite kernel from a radial basis function in such a way that shifts of the kernel can be linearly combined to yield a div-free or curl-free interpolant of the underlying field.  
This technique has the added benefits that it gives a well-posed interpolant for scattered data, is devoid of any coordinate singularities, and naturally allows a scalar potential for the field to be extracted{~\cite{FuselierWright2009}}.

When using a smooth RBF that features a shape parameter, $\ep$, to construct these div-free or curl-free kernels, one often finds that the best accuracy of the interpolated field is achieved when $\ep$ is small, corresponding to a flat kernel, but that the direct way of computing the interpolant (often called RBF-Direct) is prohibitively ill-conditioned.  This is exactly analogous to the standard RBF interpolation problem for scalar functions.  In the scalar setting, three distinctly different numerical algorithms have been presented thus far in the literature to bypass this ill-conditioning and open up the complete range of $\ep$ that can be considered.  These are the RBF-RA method~\cite{FoWr,WF17}, the RBF-QR method~\cite{FLF11,FornbergPiret2007,FaMC12}, and the RBF-GA method~\cite{FoLePo13}.  The present short note focuses on the RBF-QR method from~\cite{FornbergPiret2007}, which is specific to RBF interpolation on the sphere.  This algorithm exploits the Mercer expansion of the scalar RBFs in terms of spherical harmonics to change the interpolation basis so that the the ill-conditioning associated with small $\ep$ is analytically removed.  It also shows that scalar RBF interpolants converge to spherical harmonic interpolants in the flat limit (i.e. $\ep\rightarrow 0$).  We show how this algorithm, which we call the vector RBF-QR, can be extended to bypass the ill-conditioning associated with surface div-free RBF interpolation for small $\ep$ and demonstrate that these interpolants converge to div-free vector spherical harmonic interpolants as $\ep\rightarrow 0$.  The algorithm also extends naturally to the surface curl-free case, but the sake of brevity we leave out these details and refer the reader to~\cite{DrakeThesis}.

The rest of the paper is organized as follows:  In Section 2, we introduce some notation and background on div-free RBF interpolation and scalar/vector spherical harmonics. Section 3 contains the description of the new vector RBF-QR algorithm and Section 4 contains a numerical example illustrating the stability of the algorithm in the flat limit.  We conclude with some brief comments in Section 5.

\section{Notation and preliminaries}
\subsection{Surface div-free  vector fields}
\label{surfaceops}
Any $C^1$ tangential velocity field $\vu$ on $\mathbb{S}^2$ that is surface div-free  can be written as
\begin{align*}
\vu = \nabla^{\ast} \times (\psi \widehat{\mathbf{r}}) = -\underbrace{\hat{\mathbf{r}} \times \nabla^{\ast}}_{\rots} \psi,
\end{align*}
where $\nabla^{\ast}$ denotes the gradient in spherical coordinates, $\hat{\mathbf{r}}$ is the unit radial vector in spherical coordinate basis, and $\psi$ is some $C^2$ scalar-valued function on $\mathbb{S}^2$. The function $\psi$ is called the \textit{stream function} and is unique up to a constant for a given $\vu$~\cite[Proposition 2.1]{FuselierNarcWardWright2009}.  

The operator $\rots$, which is sometimes called the surface-curl operator or the \emph{rot} operator, can be written entirely in extrinsic (Cartesian) coordinates as $\rotc = \hat{\vn} \times \nabla$, where $\nabla$ is the standard $\mathbb{R}^3$ gradient in the Cartesian basis, and $\hat{\vn}$ is the unit normal vector to $\mathbb{S}^2$ in the Cartesian basis~\cite{NarcWardWright}.  Here we have dropped the $\ast$ from $\rotc$ to indicate it is defined in extrinsic coordinates. For scalar-valued functions $\psi$ that can be extended smoothly from $\mathbb{S}^2$ to $\mathbb{R}^3$, we can generate a tangential velocity field that is surface div-free using $\vu = \rotc \psi$.  The operator $\rotc$ can be simplified further by noting that $\hat{\vn}$ at a point $\vx = (x,y,z)$ on $\mathbb{S}^2$ is just $\vx$. In what proceeds, we use the following notation for the operator $\rotc$:
\begin{align}
\rotc_{\vx} = \underbrace{\begin{pmatrix}
		0 & -z & y\\
		z &0 & -x\\
		-y &x & 0 \end{pmatrix}}_{\ds \Qx}\nabla_{\vx},
\label{eq:rotc}
\end{align}
where $\Qx$ applied to a vector in $\mathbb{R}^3$ gives the cross product of $\vx$ with that vector and the subscripts  are used to indicate what variable each operator is applied to.  We note that one of the main benefits of working in extrinsic coordinates is that artificial coordinate singularities (e.g.\ pole singularities) can be avoided.

\subsection{Surface div-free  RBF interpolation}
\label{RBFVS}
\begin{figure}[h!]
	\centering
	\begin{tabular}{cc}
	\includegraphics[width=0.38\textwidth]{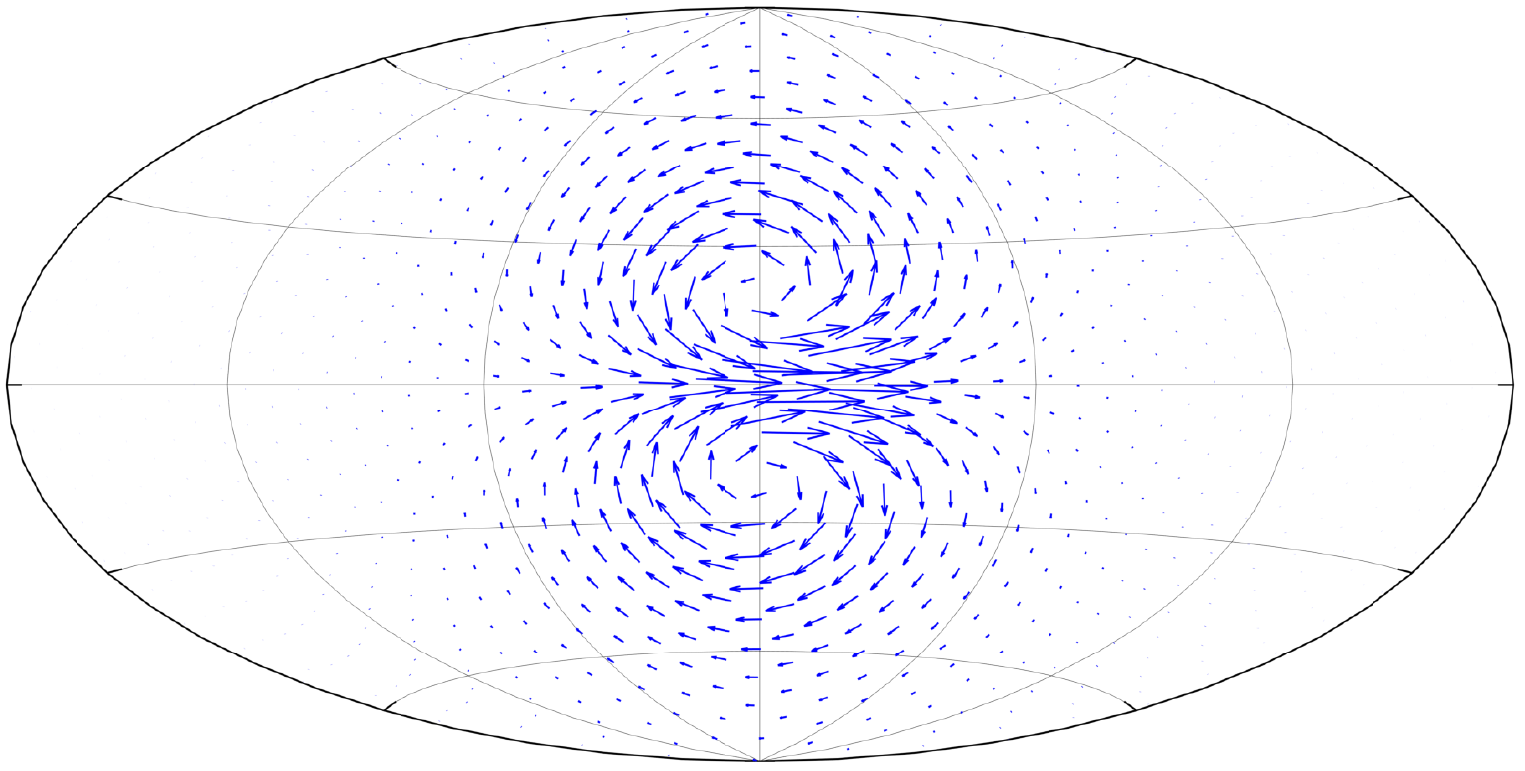} & 
	\includegraphics[width=0.38\textwidth]{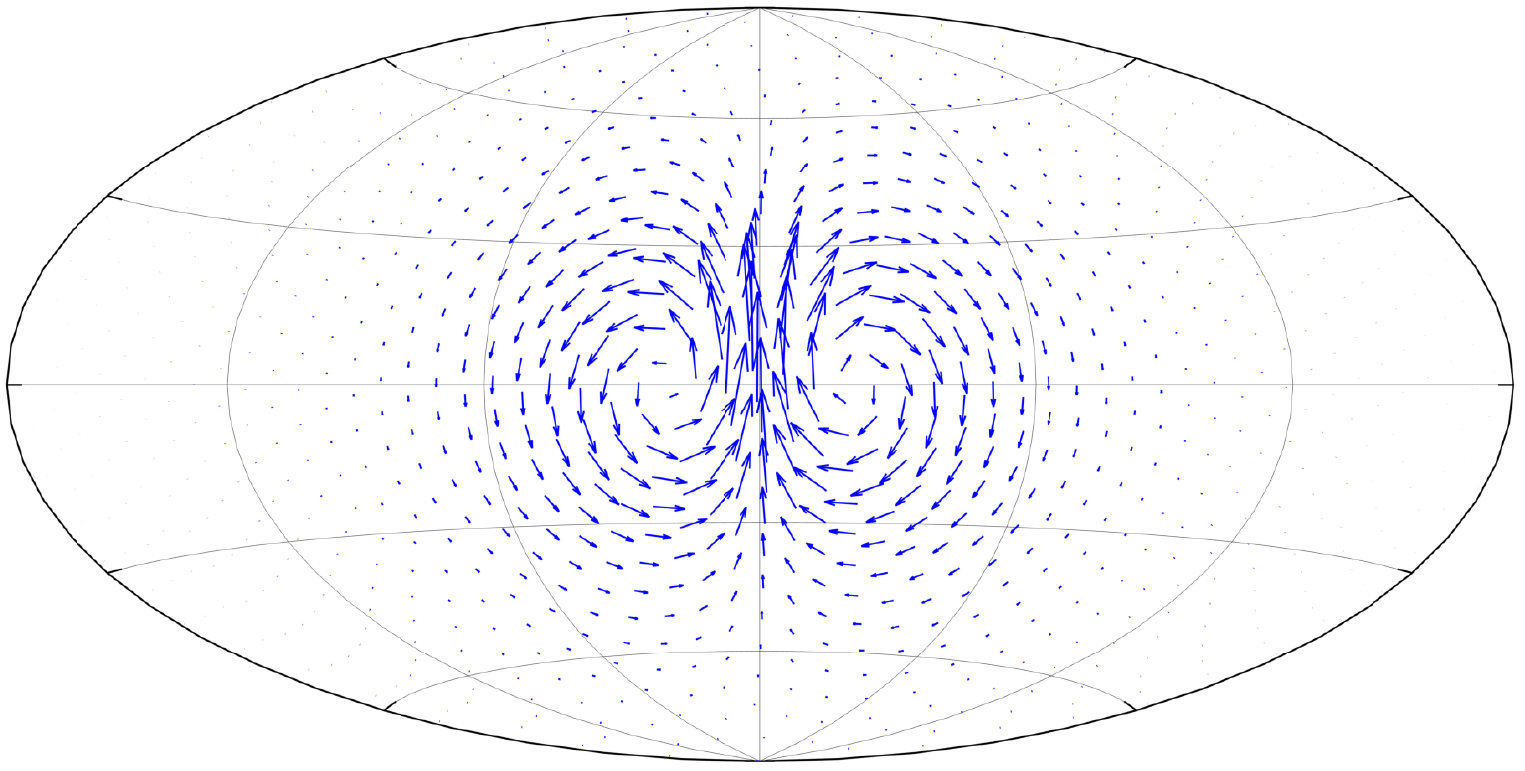} \\
	{(a)} & {(b)}
	\end{tabular}
	\caption{Illustration of the matrix-valued surface div-free  kernel \eqref{PsiDiv} using the {Hammer-Aitoff projection of the sphere~\cite{snyder1997flattening}}:  (a) $\Psidiv$ applied to the zonal unit vector at $\vy = (1,0,0)$, (a) $\Psidiv$ applied to the meridional unit vector at $\vy = (1,0,0)$.}
	\label{DivTanBas}
\end{figure}

{Surface div-free RBF interpolation is similar to scalar RBF interpolation in the sense that one constructs interpolants from linear combinations of shifts of a kernel at each of the given data sites; for a review of scalar RBF interpolation, we refer the reader to~\cite{FassBook}.  The difference between the two approaches is that surface div-free RBF interpolants use a matrix-valued kernel $\Psidiv$ whose columns are surface div-free.}  A detailed discussion of the construction of $\Psidiv$ is given in~\cite{NarcWardWright}.  For the sake of brevity, we do not repeat this derivation, but only present the final result.   Let $\phi$ be a scalar-valued, radial kernel on $\mathbb{R}^3$ with at least two continuous derivatives and let $\vx,\vy\in\mathbb{S}^2$, then $\Psidiv$ is constructed using the operator $\rotc$ in \eqref{eq:rotc} as
\begin{equation}
	\label{PsiDiv}
		\Psidiv(\vx,\vy) = \rotc_{\vx}^{\phantom{T}} \rotc_{\vy}^{T}\phi\left(\|\vx-\vy\|\right) = 
		-\Qx\left(\nabla_{\vx}^{\phantom{T}} \nabla_{\vy}^T\phi\left(\|\vx-\vy\|\right)\right)\Qy = \Qx\left(\nabla_{\vx}^{\phantom{T}} \nabla_{\vx}^T\phi\left(\|\vx-\vy\|\right)\right)\Qy,
\end{equation}
where $\| \cdot \|$ denotes the vector two-norm and we have used the fact that the $Q$ matrix in \eqref{eq:rotc} is anti-symmetric and  $\nabla_{\vy}^T\phi\left(\|\vx-\vy\|\right) = -\nabla_{\vx}^T\phi\left(\|\vx-\vy\|\right)$.  It is straightforward to show that for any $\vc\in\mathbb{R}^3$ that is tangent to $\mathbb{S}^2$ at $\vy$, the vector $\Psidiv(\vx,\vy)\vc$ is surface div-free  in $\vx$ and centered at the point $\vy$; see Figure \ref{DivTanBas} for an illustration.

Let $\{\vy_j\}_{j=1}^n$ be a distinct set of nodes on $\mathbb{S}^2$ and $\{\vu_j\}_{j=1}^n$ be samples of a surface div-free  tangent vector field sampled on these nodes.  The surface div-free  RBF interpolant to this data takes the form 
\begin{equation}
	\vs(\vx)=\sum_{j=1}^n\Psidiv(\vx,\vy_j)\vc_j,
	\label{TanDivInterp}
\end{equation}
where the coefficients $\vc_j \in \mathbb{R}^3$ are tangent to $\mathbb{S}^2$ at $\vy_j$, which is necessary to make the interpolation problem well-posed.  
A na\"ive approach to solving for the $\{\vc_j\}_{j=1}^n$ by imposing  {$\vs(\vy_j) = \vu_j$}, $j=1,\ldots,n$ will lead to a singular system of equations since each $\vc_j$ (and correspondingly $\vu_j$) can be expressed using only two degrees of freedom rather than three.  To see this, let $\{\vd_j,\ve_j,\vn_j\}$ be orthonormal vectors at the node $\vy_j$, where  $\vn_j$ is the outward normal to $\mathbb{S}^2$, $\ve_j$ is a unit tangent vector, and $\vd_j=\vn_j\times \ve_j$.  Then we can write $\vc_j$ in this basis as $\vc_j = \alpha_j\vd_j+\beta_j\ve_j$, where $\alpha_j  = \vd_j^T\vc_j$ and $\beta_j = \ve_j^T\vc_j$.  
Using this result, we can express \eqref{TanDivInterp} as 
\begin{align}
	\vs(\vx)=\sum_{j=1}^n\Psidiv(\vx,\vy_j)\underbrace{\left[\alpha_j\vd_j+\beta_j\ve_j\right]}_{\ds \mathbf{c}_j},
	\label{rbfvsinterp}
\end{align}
and write the interpolation conditions as $\vd_i^T\vs(\vy_i) = \vd_i^T\vu_i=:\gamma_i$ and  $\ve_i^T\vs(\vy_i) = \ve_i^T\vu_i=:\delta_i$, which leads to the $2n$-by-$2n$ system of equations
\begin{align}
	\sum_{j=1}^n\underbrace{\left(\begin{bmatrix}\vd_i^T\\
			\ve_i^T\end{bmatrix}\Psidiv(\vy_i,\vy_j)\begin{bmatrix}\vd_j &
			\ve_j\end{bmatrix}\right)}_{\ds A_{i,j}}\begin{bmatrix}\alpha_j\\
		\beta_j\end{bmatrix}=\begin{bmatrix}\gamma_i\\
		\delta_i\end{bmatrix}, \quad 1\leq i \leq n.
	\label{PsiDivInterp}
\end{align}
The interpolation matrix that arises from this system (with entries given by $A_{i,j}$) is positive definite if $\Psidiv$ is constructed from an appropriately chosen scalar-valued $\phi$~\cite{NarcWardWright}, such as any of those listed in Table\ref{RBFSPHEx}.  Simplifications of the entries of this matrix in terms of derivatives of $\phi$ are given in \cite{FuselierWright2009}.   

In this work, we choose $\vd_j$ and $\ve_j$, for the point $\vy_j = (x_j,y_j,z_j)$, to be the standard meridional and zonal vectors, which can be expressed in Cartesian coordinates as
\begin{equation}
	\label{tanbasvec}
	\vd_j=\frac{1}{\sqrt{1-z_j^2}}\begin{bmatrix}-z_jx_j\\
		-z_jy_j\\
		1-z_j^2 \end{bmatrix}, \quad \ve_j=\frac{1}{\sqrt{1-z_j^2}}\begin{bmatrix}-y_j\\
		x_j\\
		0 \end{bmatrix}.
\end{equation}
If $\vy_j=[0,0, \pm1]$, then we can pick any orthogonal vectors in the $xy$-plane.

We conclude by noting that once the coefficients $\vc_j$ are determined in \eqref{rbfvsinterp}, a stream-function $\psi(\vx)$ can obtained for the interpolated field using \eqref{PsiDiv}
\begin{align}
\vs(\vx)=\underbrace{\Qx\nabla_{\vx}}_{\ds \rotc_{\vx}}\biggl(\underbrace{\sum_{j=1}^n\nabla_{\vx}^T\phi\left(\|\vx-\vy_j\|\right)\Qyj\vc_j}_{\ds \psi(\vx)}\biggr).
\label{eq:rbf_sf}
\end{align}
Approximation results for $\psi$ in reconstructing the underlying stream function for the field are given in~\cite{FuselierWright2009}.

Similar to scalar RBF interpolation, we refer to the method of computing the div-free  RBF interpolant by solving the system \eqref{PsiDivInterp} as RBF-Direct.  For $\Psidiv$ that are built from scalar kernels that depend on a shape parameter $\ep$, like those in Table \ref{RBFSPHE}, the RBF-Direct approach becomes extremely ill-conditioned as $\ep\rightarrow 0$.  However, as in the scalar RBF case, the interpolant remains well-conditioned and the RBF-QR algorithm allows us to stably compute it.

\subsection{Spherical harmonics}
The RBF-QR algorithm presented in Section 3 relies on spherical harmonic expansions, both scalar and vector ones, so we briefly review these expansions here.  We refer the reader to more detailed discussions in~\cite{AtkinsonHan2012} and~\cite{Swarzt1993}, for the scalar and vector case, respectively.

Let $Y_\mu^\nu$ denote the scalar spherical harmonic of integer degree $\mu\geq 0$ and integer order $-\mu\leq \nu\leq \mu$ on $\mathbb{S}^2$. These functions are the eigenfunctions of the Laplace-Beltrami operator $\Delta_{\mathbb{S}^2}$, i.e., $\Delta_{\mathbb{S}^2}Y_\mu^\nu=-\mu(\mu+1)Y_\mu^\nu$.  We use the real-form of the spherical harmonics functions in Cartesian coordinates, which for $\vx=(x,y,z)\in \mathbb{S}^2$ are given as
\begin{align}
	Y_\mu^\nu(\vx) = \begin{cases}
		\sqrt{\frac{2\mu+1}{4\pi}}\sqrt{\frac{(\mu-\nu)!}{(\mu+\nu)!}}P_\mu^\nu(z)\cos\left(\nu \tan^{-1}\left(\frac{y}{x}\right)\right),& \text{ } \nu=0,1,\dots,\mu, \\ 
		\sqrt{\frac{2\mu+1}{4\pi}}\sqrt{\frac{(\mu-\nu)!}{(\mu+\nu)!}}P_\mu^\nu(z)\sin\left(-\nu \tan^{-1}\left(\frac{y}{x}\right)\right),& \text{ } \nu= -\mu,\dots,-1
	\end{cases},
	\label{ssph}
\end{align}
where $P_\mu^\nu(z)$ are the associated Legendre functions of degree $\mu$ and order $\nu$. The spherical harmonics form a complete, orthonormal set of basis functions for the space of square-integrable functions on $\mathbb{S}^2$, which we denote by $L^2(\mathbb{S}^2)$~\cite{AtkinsonHan2012}. Thus, any function $f\in L^2(\mathbb{S}^2)$ can be uniquely represented as
\begin{align}
	f(\vx) = \sum_{\mu=0}^\infty \sum_{\nu=-\mu}^\mu \hat{f}_{\mu,\nu}Y_\mu^\nu (\vx),\; \text{where}\;\;\hat{f}_{\mu,\nu}= \int_{\mathbb{S}^2}f(\vx) Y_{\mu}^{\nu}(\vx)dS.
	\label{SPHE}
\end{align} 

Vector spherical harmonics are the vectorial analogue of scalar spherical harmonics and can be used for representing vector-valued functions on $\mathbb{S}^2$. There are three types of these vector harmonic functions: one that is normal to $\mathbb{S}^2$, one that is tangent to $\mathbb{S}^2$ and surface curl-free, and one that is tangent to $\mathbb{S}^2$ and surface div-free~\cite{Swarzt1993}. We focus on the latter as they are the ones used in this paper.  Both tangent vector spherical harmonics are the eigenfunctions of the vector Laplace-Beltrami operator.

The surface div-free  vector spherical harmonics can be constructed in Cartesian coordinates by applying operator $\rotc$ in \eqref{eq:rotc} to the scalar spherical harmonic functions \eqref{ssph} as follows~\cite{Swarzt1993}:
\begin{equation}
	\vw_\mu^\nu(\vx)=\rotc_{\vx} Y_\mu^\nu(\vx),\; \mu>0,\; -\mu \leq \nu \leq \mu.
	\label{divVSPH}
\end{equation}
These can be normalized as $\overline{\vw}_\mu^\nu = \frac{1}{\mu(\mu+1)}\vw_\mu^\nu(vx)$ so they are orthonormal with respect to the $L^2(\mathbb{S}^2)$-vector inner product $\langle \mathbf{f},\mathbf{g}\rangle = \int_{\mathbb{S}^2}\mathbf{f}^T\mathbf{g}dS$,
where $\mathbf{f}$ and $\mathbf{g}$ are {tangent} vector fields on $\mathbb{S}^2$.  
The  surface div-free  vector spherical harmonics form a complete orthonormal set of vector basis functions for the space of surface div-free  vector fields that are tangent to $\mathbb{S}^2$ and square integrable.  We denote this space as $L_{\rm div}(\mathbb{S}^2)$. Any function $\vu\in L_{\rm div}^2(\mathbb{S}^2)$ can be expanded in these harmonics as 
\begin{equation}
	\vu(\vx)=\sum_{\mu=1}^\infty\sum_{\nu=-\mu}^\mu\hat{u}_{\mu,\nu}\overline{\vw}_\mu^\nu(\vx),\; \text{where}\;\; \hat{u}_{\mu,\nu}=\int_{\mathbb{S}^2}\mathbf{f(\vx)}^T\overline{\vw}_\mu^\nu(\vx)dS.
	\label{eq:vshpex}
\end{equation}
Note here that the sum in this expansion excludes $\mu=0$ since $Y_0^0$ is annihilated by $\rotc$.  

As discussed in Section~\ref{RBFVS}, the interpolation problem on the sphere requires that we utilize the tangent basis vectors at $\vx\in\mathbb{S}^2$, i.e., $\vd_\vx$ and $\ve_\vx$. Therefore, it is relevant to introduce the following notation for the non-normalized surface div-free  and curl-free vector spherical harmonics in terms of these basis vectors:
\begin{align}
	\textnormal{meridional, div-free :}\; G_\mu^\nu(\vx)=\vd_\vx^T  \rotc_{\vx} Y_\mu^\nu(\vx), \qquad  
	\textnormal{zonal, div-free :}\; H_\mu^\nu(\vx)=\ve_\vx^T  \rotc_{\vx} Y_\mu^\nu(\vx). 
	\label{eq:vsph_components}
\end{align}


\subsection{Spherical harmonic expansions of smooth radial kernels}
The central idea of the scalar RBF-QR method is to replace the standard basis consisting of shifts of radial kernels by an equivalent, but much better conditioned basis that spans the same space in the case of small $\ep$.  For the sphere, this is done by exploiting properties of the Mercer expansion of smooth radial kernels on the sphere in terms of spherical harmonics.  As discussed in~\cite{FornbergPiret2007}, these expansions are given by
\begin{equation}
	\phi(\|\vx-\vy_j\|)=\sum_{\mu=0}^\infty \sideset{}{'}\sum_{\nu=-\mu}^{\mu}\{c_{\mu,\ep}\ep^{2\mu}Y_\mu^\nu(\vy_j)\}Y_\mu^\nu(\vx),
	\label{RBFSPHEx}
\end{equation}
where the symbol $\sum'$ denotes that the $\nu=0$ term is halved.  Note that the spherical harmonic coefficients are independent of $\nu$, which follows from the Funk-Hecke formula~\cite{Muller1966,AtkinsonHan2012} for zonal functions. 
Table~\ref{RBFSPHE} lists the coefficients $c_{\mu,\ep}$ for many common radial kernels. These were first computed by Hubbert and Baxter~\cite{HubBax} for all the radial kernels listed in this table except for the IQ, which was given in~\cite{FornbergPiret2007}.

\begin{table}[h!]
	\centering
	\begin{tabular}{|c|c|c|}
		\hline
		\textbf{Radial kernel} & \textbf{Expression} & \textbf{Expansion coefficient, $c_{\mu,\ep}$}                                                                            \\ \hline
		Multiquadric (MQ)                                             & $\phi(r) = (1 + (\ep r)^2)^{\frac12}$ & $\frac{-2\pi(2\ep^2+1+(\mu+1/2)\sqrt{1+4\ep^2})}{(\mu+3/2)(\mu+1/2)(\mu-1/2)}\left(\frac{2}{1+\sqrt{4\ep^2+1}}\right)^{2\mu+1}$ \\ \hline
		Inverse multiquadric (IMQ)                                            & $\phi(r) = (1 + (\ep r)^2)^{-\frac12}$ & $\frac{4\pi}{(\mu+1/2)}\left(\frac{2}{1+\sqrt{4\ep^2+1}}\right)^{2\mu+1}$                                                                                                                    \\ \hline
		Inverse quadratic (IQ)                                             & $\phi(r) = (1 + (\ep r)^2)^{-1}$ & $\frac{4\pi^{3/2}\mu!}{\Gamma(\mu+\frac{3}{2})(1+4\ep^2)^{\mu+1}} \prescript{}{2}{F}_1(\mu+1,\mu+1;2\mu+2;\frac{4\ep^2}{1+4\ep^2})$              \\ \hline
		Gaussian (GA)                                             & $\phi(r) = \exp(-(\ep r)^2)$ & $\frac{4\pi^{3/2}}{\ep^{2\mu+1}}e^{-2\ep^2}I_{\mu+1/2}(2\ep^2)$                                                                 \\ \hline
	\end{tabular}
	\caption{\label{RBFSPHE} Coefficients in the spherical harmonic expansion \eqref{RBFSPHEx} for various smooth radial kernels on the sphere. For the IQ kernel, $\prescript{}{2}{F}_1(\dots)$ denotes the hypergeometric function, and for the GA kernel, $I_{\mu+1/2}$ denotes the Bessel function of the second kind.}
\end{table}
It is important to note that the coefficients listed in Table~\ref{RBFSPHE} can be calculated without the loss of any significant digits caused by numerical cancellations, even as $\ep\rightarrow 0$ ({see~\cite{FornbergPiret2007} for a discussion}).

The central idea behind the vector RBF-QR algorithm will also be to replace the matrix-valued basis with a better basis built from vector spherical harmonic expansions.  We can use \eqref{PsiDiv} and \eqref{RBFSPHEx} to expand the matrix-valued kernel $\Psidiv$ in terms of vector spherical harmonics as follows:
\begin{align}
		\Psidiv(\vx,\vy_j)
		=\sum_{\mu=1}^\infty \sideset{}{'}\sum_{\nu=-\mu}^\mu \{\ep^{2\mu}c_{\mu,\ep}\rotc_{\vx} Y_\mu^\nu(\vx)\}\left(\left.\rotc_{\vy}Y_\mu^\nu(\vy)\right|_{\vy=\vy_j}\right)^T=\sum_{\mu=1}^\infty \sideset{}{'}\sum_{\nu=-\mu}^\mu \{\ep^{2\mu}c_{\mu,\ep}\vw_\mu^\nu(\vx)\}\left(\vw_\mu^\nu(\vy_j)\right)^T.
		\label{PsiDivVSPH}
\end{align}
Here we have used the non-normalized div-free  vector spherical harmonics defined in~\eqref{divVSPH} for simplicity.

The separation of the expansion \eqref{PsiDivVSPH} of a smooth matrix-valued kernel $\Psidiv$ in terms of increasing powers of $\ep$ is essential to the RBF-QR algorithm.  It is these powers of $\ep$ and not the coefficients $c_{\mu,\ep}$ that lead ill-conditioning in the matrix-valued basis for small $\ep$ in the RBF-Direct method.  The RBF-QR algorithm analytically factors out the effects of these powers of $\ep$ from the basis.

\section{Vector RBF-QR algorithm}\label{ch:RBFQRV}

Recall from Section~\ref{RBFVS} that in order to interpolate div-free  vector fields tangent to the sphere with the matrix-valued div-free  interpolant, we must represent the coefficient vectors and target field samples in terms of the orthonormal tangent basis vectors (e.g., using~\eqref{tanbasvec}). In~\eqref{PsiDivInterp} we saw that this is equivalent to representing the kernel $\Psidiv$ in terms of these basis vectors. We will denote this kernel as $\Psidivt$:
\begin{equation}
	\label{PsiDivTan}
	\Psidivt(\vx,\vy_j)=
	\begin{bmatrix}\vd_\vx^T\\
		\ve_\vx^T\end{bmatrix}\Psidiv(\vx,\vy_j)\begin{bmatrix}\vd_j & \ve_j
	\end{bmatrix}.
\end{equation}
Using~\eqref{PsiDivVSPH} on the right-hand side of~\eqref{PsiDivTan} gives the expansion
\begin{equation}
\Psidivt(\vx,\vy_j)= \sum_{\mu=1}^\infty \sideset{}{'}\sum_{\nu=-\mu}^\mu \left(\ep^{2\mu}c_{\mu,\ep}\begin{bmatrix}
		\vd_\vx^T\\
		\ve_\vx^T
	\end{bmatrix}\vw_\mu^\nu(\vx)\right)\left(\vw_\mu^\nu(\vy_j)\right)^T\begin{bmatrix}\vd_j & \ve_j\end{bmatrix}.\end{equation}
This is a $2$-by-$2$ matrix whose entries are in terms of the meridional and zonal div-free  vector spherical harmonics \eqref{eq:vsph_components}:
\begin{equation}	
	\Psidivt(\vx,\vy_j)=
	\begin{bmatrix}
	\begin{array}{c|c}
			(a) &(b) \\
			\hline
			(c) &(d) \\
	\end{array}
	\end{bmatrix}, \; \textnormal{where}\label{1PsiDivMat}
\end{equation}
\begin{align*}
	\text{(a)}\; \sum_{\mu=1}^\infty \sideset{}{'}\sum_{\nu=-\mu}^\mu \{\ep^{2\mu}c_{\mu,\ep}G_\mu^\nu(\vx)\} G_\mu^\nu(\vy_j) &\qquad 
	\text{(b)}\; \sum_{\mu=1}^\infty \sideset{}{'}\sum_{\nu=-\mu}^\mu \{\ep^{2\mu}c_{\mu,\ep}G_\mu^\nu(\vx)\} H_\mu^\nu(\vy_j)\\
	\text{(c)}\; \sum_{\mu=1}^\infty \sideset{}{'}\sum_{\nu=-\mu}^\mu \{\ep^{2\mu}c_{\mu,\ep}H_\mu^\nu(\vx)\} G_\mu^\nu(\vy_j) &\qquad  
	\text{(d)}\; \sum_{\mu=1}^\infty \sideset{}{'}\sum_{\nu=-\mu}^\mu \{\ep^{2\mu}c_{\mu,\ep}H_\mu^\nu(\vx)\} H_\mu^\nu(\vy_j).
\end{align*} 

The goal of the vector RBF-QR algorithm is to express the space spanned by the columns of the $2$-by-$2n$ array containing the shifts of the div-free  matrix valued kernel,
$
\begin{bmatrix}
\Psidivt(\vx,\vy_1) & \cdots & \Psidivt(\vx,\vy_n)
\end{bmatrix},
$
using a basis that has the ill-conditioning associated with small $\ep$ removed.  To this end, we first use \eqref{1PsiDivMat} to write this array (now in transposed form) as 
\begin{equation} \label{DivMat}
\begin{bmatrix}
		\Psidivt(\vx,\vy_1)^{T}\\ 
		\vdots\\ 
		\Psidivt(\vx,\vy_n)^{T}
	\end{bmatrix}=
\begin{bmatrix}
		\begin{bmatrix}\begin{array}{c|c}
				\ds\sum_{\mu=1}^\infty \sideset{}{'}\sum_{\nu=-\mu}^\mu \{\ep^{2\mu}c_{\mu,\ep}G_\mu^\nu(\vx)\} G_\mu^\nu(\vy_1) & \ds\sum_{\mu=1}^\infty \sideset{}{'}\sum_{\nu=-\mu}^\mu \{\ep^{2\mu}c_{\mu,\ep}H_\mu^\nu(\vx)\} G_\mu^\nu(\vy_1)\\
				\hline
				\ds\sum_{\mu=1}^\infty \sideset{}{'}\sum_{\nu=-\mu}^\mu \{\ep^{2\mu}c_{\mu,\ep}G_\mu^\nu(\vx)\} H_\mu^\nu(\vy_1) &\ds\sum_{\mu=1}^\infty \sideset{}{'}\sum_{\nu=-\mu}^\mu \{\ep^{2\mu}c_{\mu,\ep}H_\mu^\nu(\vx)\} H_\mu^\nu(\vy_1) \\
			\end{array}\end{bmatrix}\\
			\vdots\\
			\begin{bmatrix}\begin{array}{c|c}
					\ds\sum_{\mu=1}^\infty \sideset{}{'}\sum_{\nu=-\mu}^\mu \{\ep^{2\mu}c_{\mu,\ep}G_\mu^\nu(\vx)\} G_\mu^\nu(\vy_n) & \ds\sum_{\mu=1}^\infty \sideset{}{'}\sum_{\nu=-\mu}^\mu \{\ep^{2\mu}c_{\mu,\ep}H_\mu^\nu(\vx)\} G_\mu^\nu(\vy_n)\\
					\hline
					\ds\sum_{\mu=1}^\infty \sideset{}{'}\sum_{\nu=-\mu}^\mu \{\ep^{2\mu}c_{\mu,\ep}G_\mu^\nu(\vx)\} H_\mu^\nu(\vy_n) &\ds\sum_{\mu=1}^\infty \sideset{}{'}\sum_{\nu=-\mu}^\mu \{\ep^{2\mu}c_{\mu,\ep}H_\mu^\nu(\vx)\} H_\mu^\nu(\vy_n) \\
				\end{array}\end{bmatrix}
			\end{bmatrix}.\end{equation}
We then rewrite this as the following infinite block matrix-matrix product,
\begin{align}
		\underbrace{
		\begin{bmatrix}
				\ds	\Psidivt(\vx,\vy_1)^T\\ 
				\vdots\\ 
				\ds\Psidivt(\vx,\vy_n)^T
		\end{bmatrix}}_{\ds \Pdiv}
		=
		\underbrace{
		\begin{bmatrix}
				c_{1,\ep}G_1^{-1}(\vy_1)  &\frac{c_{1,\ep}}{2}G_1^0(\vy_1)  &c_{1,\ep} G_1^1(\vy_1) \cdots\\
				c_{1,\ep} H_1^{-1}(\vy_1)  &\frac{c_{1,\ep}}{2}H_1^0(\vy_1)  &c_{1,\ep} H_1^1(\vy_1) \cdots\\
				\vdots  &\vdots  &\vdots \\
				c_{1,\ep} G_1^{-1}(\vy_n)  &\frac{c_{1,\ep}}{2}G_1^0(\vy_n)  &c_{1,\ep} G_1^1(\vy_n) \cdots\\
				c_{1,\ep} H_1^{-1}(\vy_n)  &\frac{c_{1,\ep}}{2}H_1^0(\vy_n)  &c_{1,\ep} H_1^1(\vy_n) \cdots\\
		\end{bmatrix}}_{\ds B^\infty}
		\underbrace{
		\begin{bmatrix}
					\ep^2  &  & \\ 
					&\ep^2  & \\ 
					&  &\ep^2\\
					&   &  &\ddots \vspace{-0.1in}\\
					& & & &\hspace{-0.1in}\ddots
		\end{bmatrix}}_{\ds E^\infty}
		\underbrace{
		\begin{bmatrix}
				G_1^{-1}(\vx) &H_1^{-1}(\vx)\vspace{1.3ex}\\ 
				G_1^0(\vx) &H_1^0(\vx)\vspace{1.3ex}\\ 
				G_1^1(\vx) &H_1^1(\vx)\\
				\vdots &\vdots\\
		\end{bmatrix}}_{\ds Y^\infty}.
		\label{PsiDivFullMat}
\end{align}

The next step is to truncate these infinite matrices to a vector spherical harmonic degree value $\mu=\mutrunc$.  There are two stipulations for this truncation degree. First, $\mutrunc$ must be large enough that the entries in $\Pdiv$ are approximated to machine precision; we refer to this value as $\mu_{\text{eps}}$. The second is that $\mutrunc$ must be large enough that the number of columns in the truncated $B^{\infty}$ matrix is greater than or equal to $2n$; we refer to this value as $\mu_0$.  This condition on $\mu_0$ is given explicitly as $\mu_0= \left\lceil \sqrt{2n+1}-1\right\rceil$ since the number of vector spherical harmonic terms in a truncated expansion of $\mu_0$ is $\mu_0(\mu_0 + 2)$.  Putting the truncation requirements together gives the  condition $\mutrunc\geq\max(\mu_0,\mu_{\text{eps}})$. We denote the truncated matrix-matrix product~\eqref{PsiDivFullMat} as 
\begin{align}
    \Pdiv \approx
    \begin{bmatrix}
		\\
		& & & B   & & &\\ 
		\\
    \end{bmatrix}
    \begin{bmatrix}
		\\
		& &  E   & &  \\ 
		\\
    \end{bmatrix}  
    \begin{bmatrix}
		\\ 
		Y   \\ 
		\\ 
    \end{bmatrix}.\label{RBFVSPHTrunc}
\end{align}
The matrix $B$ has $m=\mutrunc(\mutrunc+2)$ columns and has the block form
\begin{align}
	B = \left[
		\begin{array}{c|c|c|c|c|c|c}
			B_1 & B_2 & \cdots & B_{\mu_0} & B_{\mu_0 + 1} & \cdots & B_{\mutrunc}
		\end{array}
		\right]
	\label{eq:BTrunc}
\end{align}
where $B_\mu$, $1\leq \mu \leq \mutrunc$, are the block matrices of size $2n$-by-$(2\mu+1)$ with block entries
\begin{equation*}
	\left(B_\mu\right)_{i,j} =
	\left\{\begin{matrix}
		\begin{bmatrix}
			c_{\mu,\ep} G_{\mu}^{j-(\mu+1)}(\vy_i)\;\\ 
			c_{\mu,\ep} H_{\mu}^{j-(\mu+1)}(\vy_i)
		\end{bmatrix}j\neq \mu+1,\vspace{2ex}\\ 
		\begin{bmatrix}
			\frac{c_{\mu,\ep}}{2} G_{\mu}^{0}(\vy_i)\;\\ 
			\frac{c_{\mu,\ep}}{2} H_{\mu}^{0}(\vy_i)
		\end{bmatrix}\; \qquad j = \mu+1,
	\end{matrix}\right. \;\;
	j=1,\ldots,2\mu+1, \;i=1,\ldots,n.
\end{equation*}
The diagonal $E$ matrix in \eqref{RBFVSPHTrunc} can be written as two square, diagonal blocks, $E = \text{diag}([E_1;E_2])$, where
\begin{align}
	E_1 = 
	\begin{bmatrix}
		& \ep^2 I_{3} \\
		& & \ep^4 I_{5} \\
		& & & \ddots \\
		& & & & \ep^{2\mu_0} I_{2\mu_0+1}
	\end{bmatrix}\; \text{and}\;
	E_2  = 
	\begin{bmatrix}
		\ep^{2\mu_0 + 2} I_{2\mu_0 + 3} \\
		& \ep^{2\mu_0 + 4} I_{2\mu_0 + 5} \\
		& & \ddots \\
		& & & \ep^{2\mutrunc} I_{2\mutrunc+1}
	\end{bmatrix} \label{VEMatrix},
\end{align}
and $I_\tau$ is the identity matrix of size $\tau$-by-$\tau$. 
Finally, the $Y$ matrix in \eqref{RBFVSPHTrunc}  is given in block form by
\begin{align*}
	Y^{T} = \left[
		\begin{array}{c|c|c|c|c|c|c}
			Y_1^{T} & Y_2^{T} & \cdots & Y_{\mu_0}^{T} & Y_{\mu_0 + 1}^{T} & \cdots & Y_{\mutrunc}^{T}
		\end{array}
		\right],
\end{align*}
where block $Y_{\mu}$ is given by
\begin{align*}
	\left(Y_\mu\right)_{j,1} = 	G_{\mu}^{j-(\mu+1)}(\vx), \quad \left(Y_\mu\right)_{j,2} = H_{\mu}^{j-(\mu+1)}(\vx),
	\quad
	j=1,\ldots,2\mu+1.
\end{align*}

In the flat limit, {the truncated basis} is still highly ill-conditioned because of the powers of $\ep$ in~\eqref{PsiDivFullMat} (recall that the expansion coefficients $c_{\mu,\ep}$ do not affect the conditioning).  However, all of these powers of $\ep$ are confined to the $E$ matrix. To develop a better conditioned basis, we need to factor out the ill-effects of the powers of $\ep$.  To simplify the description of this step of the algorithm, we make the assumption that $n = \mu_0(\mu_0 + 2)/2$ so that the sub-matrix block 
\begin{align}
\left[\begin{array}{c|c|c|c}
	 B_1 & B_2 & \cdots & B_{\mu_0} 
	\end{array}
\right]
\label{eq:Bsubblock}
\end{align}
of $B$ in \eqref{eq:BTrunc} and diagonal matrix $E_1$ in \eqref{VEMatrix} are square and of size $2n$-by-$2n$.  We also assume that the interpolation node set $\{\vy_j\}_{j=1}^n$ is unisolvent with respect to the vector spherical harmonics of degree $\mu_0$ so that \eqref{eq:Bsubblock} is invertible.  {This restriction means that the interpolation matrix associated with the vector spherical harmonic basis of degree $\mu_0$ is invertible.  When using ``scattered'' node sets, we have never encountered a situation where the point set fails to be unisolvent.}
 
The final step of the vector RBF-QR algorithm starts with a QR factorization on $B$ in \eqref{RBFVSPHTrunc}, which gives
\begin{equation}
	\widetilde{P}_{\rm div} \approx Q\underbrace{\left[R_1\,|\,R_2\right]}_{\displaystyle R}\left[
	\begin{array}{c|c}
		E_1 & \\
		\hline
		& E_2 \\
	\end{array}
	\right]Y.
	\label{QRV1}
\end{equation}
Here we have partitioned $R$ into $R_1$ and $R_2$, where $R_1$ is $2n$-by-$2n$ upper-triangular, and $R_2$ is an $2n$-by-$(m-2n)$ full matrix.  Since the sub-block \eqref{eq:Bsubblock} is invertible by our assumption on the node set, we know $R_1$ is also invertible.  This together with the block structure of $E$ allows us to re-write \eqref{QRV1} as 
\begin{align}
	\widetilde{P}_{\rm div} \approx Q R_1 \left[I_{2n}\,|\,R_1^{-1}R_2\right]\left[
	\begin{array}{c|c}
		E_1 & \\
		\hline
		& E_2 \\
	\end{array}
	\right]Y = 
	Q R_1 \left[E_1\,|\,R_1^{-1}R_2E_2\right]Y = Q R_1 E_1 \underbrace{\left[I_{2n}\,|\,E_1^{-1}R_1^{-1}R_2E_2\right]}_{\displaystyle \widetilde{B}_{\rm div}}Y.
	\label{VQR3} 
\end{align}
Since left matrix-multiplication is just a linear combination of the columns of the matrices, it follows from this new expression that any element in the span of the columns of $\widetilde{P}_{\rm div}^T$ (i.e.,\ the span of the original basis containing shifts of $\Psidivt$) can be represented to machine precision by a linear combination of the columns of $(\widetilde{B}_{\rm div}Y)^T$. 

The form of $\widetilde{B}_{\rm div}Y$ in \eqref{VQR3} is still not directly amenable to computations with small $\ep$ because it involves computing $E_1^{-1}$.  However, this matrix and $E_2$ are diagonal so we can analytically remove the division by small $\ep$ using~\cite[Lemma 5.1.2]{HornJohnson1991} to arrive at
\begin{equation}
	\widetilde{B}_{\rm div}Y = \Bigl[I_{2n}\,|\,E_1^{-1}R_1^{-1}R_2E_2\Bigr] = \Bigl[I_{2n}\,|\,\left(R_1^{-1} R_2\right)\circ \underbrace{(E_1^{-1}J_{2n,m-2n}E_2)}_{\displaystyle \widetilde{E}}\Bigr],
	\label{FinalVGoodBasis}
\end{equation}
where $J_{\sigma,\tau}$ is a $\sigma$-by-$\tau$ matrix of $1$'s, $\circ$ denotes the Hadamard product (or entry-wise multiplication), and the entries of $\widetilde{E}$ can be determine explicitly as
\begin{equation} 
	\label{VEtild}
\widetilde{E} = 	\left[
	\begin{array}{c|c|cc|c}
		\ep^{2\mu_0} J_{3,2\mu_0 + 3} & \ep^{2\mu_0 + 2} J_{3,2\mu_0 + 5} & \cdots & \cdots & \ep^{2\mutrunc-2} J_{3,2\mutrunc + 1} \\
		\hline
		\ep^{2\mu_0-2} J_{5,2\mu_0 + 3} & \ep^{2\mu_0} J_{5,2\mu_0 + 5} & \cdots & \cdots & \ep^{2\mutrunc-4} J_{5,2\mutrunc + 1} \\
		\hline
		\vdots & \vdots & \ddots & \ddots & \vdots \\
		\hline
		\ep^{4} J_{2\mu_0 - 1,2\mu_0 + 3} & \ep^{6} J_{2\mu_0 - 1,2\mu_0 + 5} & \cdots & \cdots & \ep^{2\mutrunc-2\mu_0+2} J_{2\mu_0 - 1,2\mutrunc + 1} \\
		\hline
		\ep^{2} J_{2\mu_0 + 1,2\mu_0 + 3} & \ep^{4} J_{2\mu_0 + 1,2\mu_0 + 5} & \cdots & \cdots & \ep^{2\mutrunc-2\mu_0} J_{2\mu_0 + 1,2\mutrunc + 1}
	\end{array}
	\right].
\end{equation}
The columns $(\widetilde{B}_{\rm div}Y)^T$ in \eqref{FinalVGoodBasis} can now be used as a stable basis for the space spanned by $\{\widetilde{{\Psi}}_{\rm div}(\cdot,\vy_j)\}_{j=1}^n$ for small $\ep$.   We note that each element of this basis consists of a div-free vector spherical harmonic of some degree $\leq \mu_0$ plus some $\mathcal{O}(\ep^2)$ combination of div-free vector spherical harmonics of degree $> \mu_0$. This implies that matrix-valued div-free  RBF interpolant will converge to a div-free vector spherical harmonic interpolant of degree $\mu_0$ in limit $\ep\rightarrow 0$, which is analogous to a scalar RBF interpolant converging to a spherical harmonic interpolant in the flat limit~\cite{FornbergPiret2007}.

%

\begin{remark}
Note that just as with the scalar RBF-QR algorithm~\cite{FornbergPiret2007}, it is possible to include the spherical harmonic coefficients $c_{\mu,\ep}$ in the diagonal matrices $E_1$ and $E_2$ and generate a similar analytical simplification for $\widetilde{E}$. This has the added advantage of removing all $\ep$ dependence in the actual QR numerical computation and allows the central part of the vector RBF-QR to be performed independent of the radial kernel used.
\end{remark}

\begin{remark}
When the condition $2n = \mu_0(\mu_0 + 2)$ is not met, the $R_1$ matrix in \eqref{QRV1} will not be square.  The algorithm then needs to be modified to move columns from $R_1$ to $R_2$ to make $R_1$ square.  This requires a similar move of the corresponding diagonals of $E_1$ to $E_2$.  In this case, the form of $\widetilde{E}$ in \eqref{VEtild} also needs to change to include some $\ep^0$ terms.
\end{remark}

\begin{remark}
{More sophisticated techniques for selecting the truncation level $\mu_{\rm trunc}$ are discussed in~\cite{KormannLasserYurova} for the scalar case and may also be adopted here.  These methods can reduce the overall computational time without reducing the approximation properties of the interpolants based on the stable basis.}
\end{remark}

\vspace{-0.1in}
\section{Numerical Example}\label{sec:NumRes}
\begin{figure}[htb]
\centering
\begin{tabular}{cc}
\includegraphics[width=0.38\textwidth]{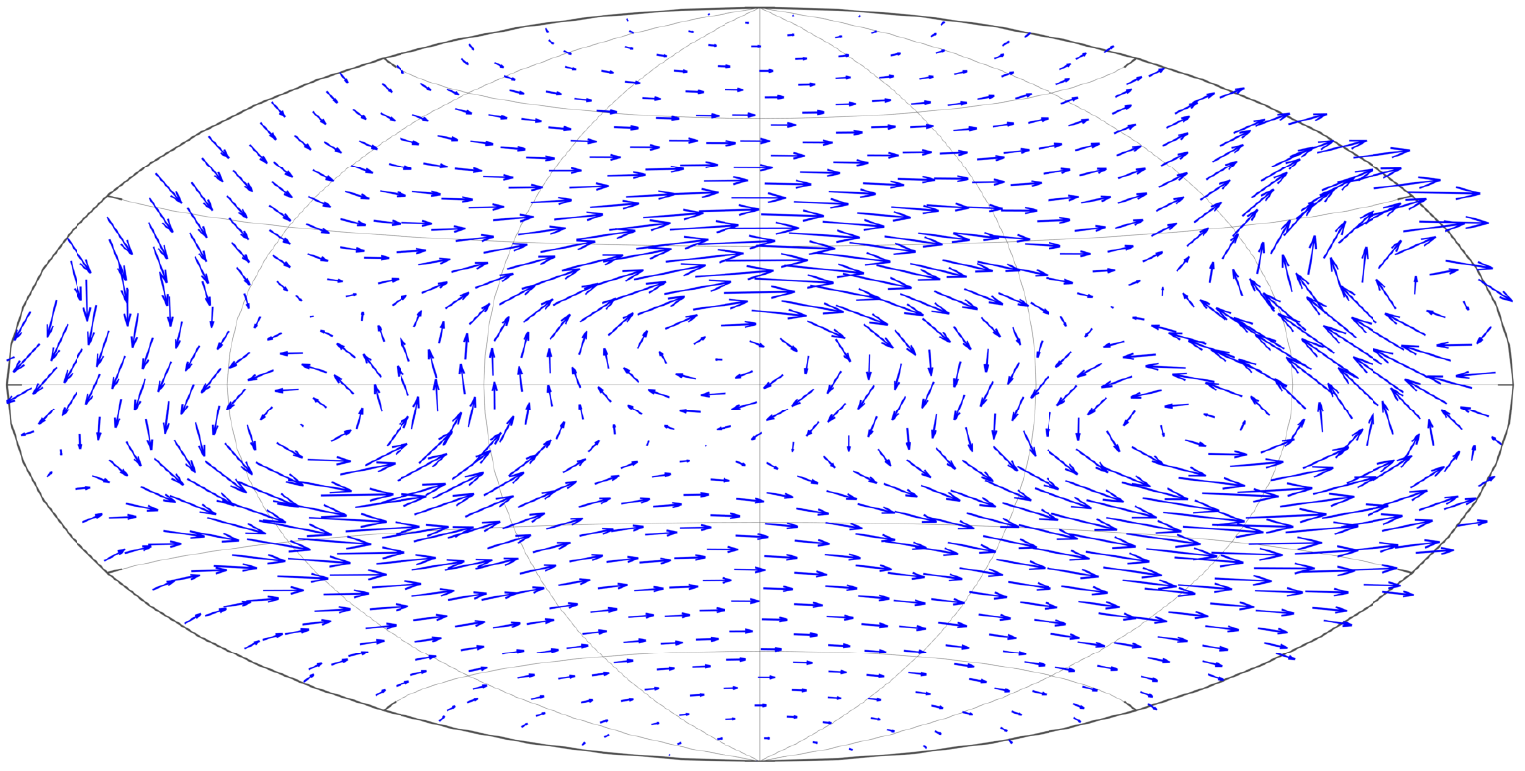} &
\includegraphics[width=0.38\textwidth]{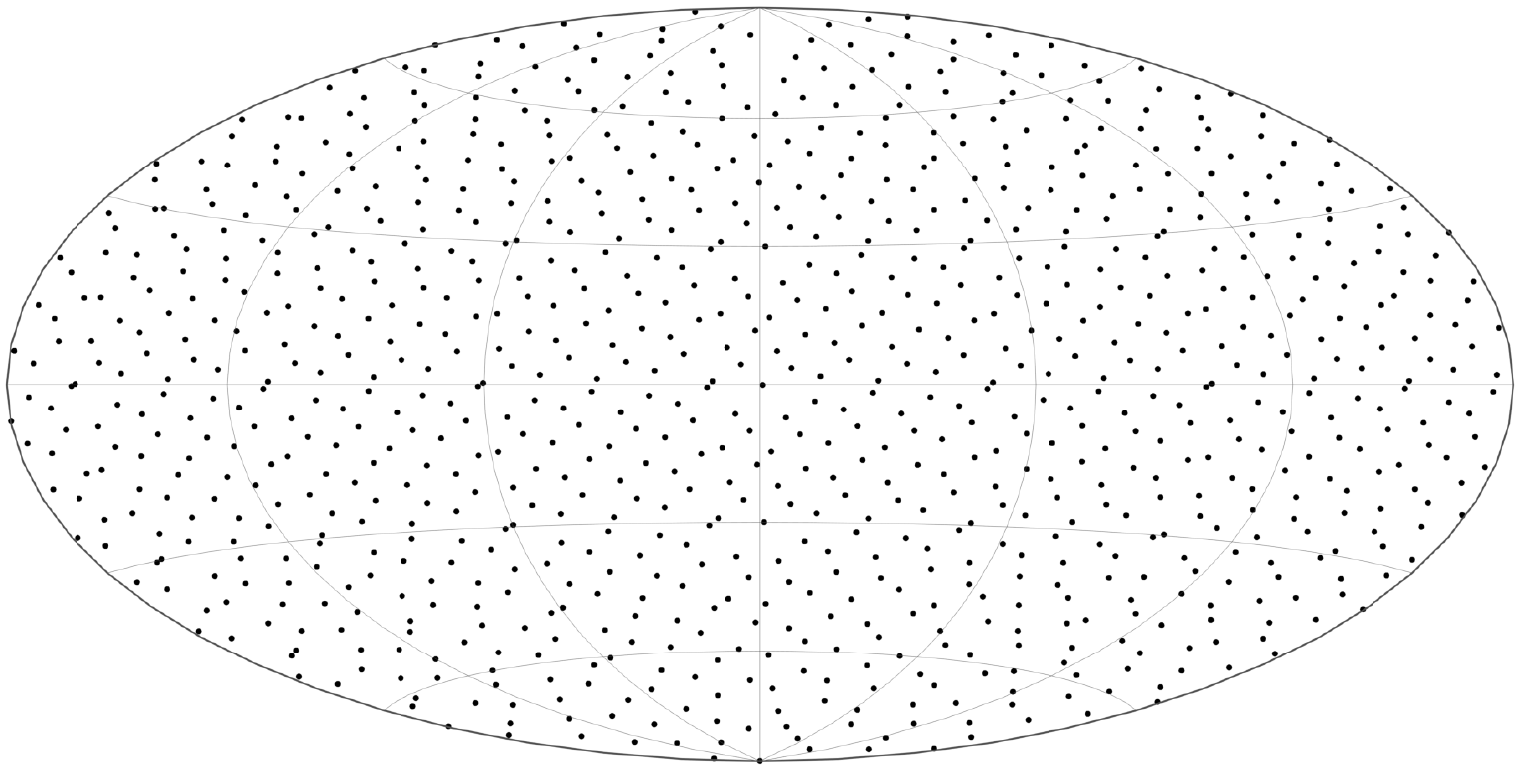} \\
{(a) Target field} & (b) Interpolation nodes
\end{tabular}
\caption{{Hammer-Aitoff projection~\cite{snyder1997flattening}} of the (a) target div-free  field and (b) $n=924$ node set used in the numerical example. \label{fig:target_field}}
\end{figure}

To test the vector RBF-QR algorithm, we use the target div-free  field displayed in Figure \ref{fig:target_field} (a).  This field is generated from the stream function 
\begin{align}
\begin{gathered}
\psi(\vx) = -3z + 2e^{-1.5((x-0.9)^2 + (y+0.1)^2) - 8(z-0.2)^2} + 3 e^{-2((x+0.7)^2 + (y-0.2)^2) - 8 (z-0.25)^2} - \\
 2.5 e^{-1.1((x+0.2)^2 + (y-0.8)^2)-8(z+0.19)^2} -  2 e^{-2.2((x+0.2)^2 + (y+1)^2) - 8(z+0.21^2)}
\end{gathered}
\label{eq:target_sf}
\end{align}
according to $\vu = \rotc_{\vx} \psi(\vx)$.  For the interpolation node set, we use the $n=924$ scattered nodes displayed in Figure \ref{fig:target_field} (b).  These are an example of the Hammersley nodes, which give well-distributed, but random  sampling points for the sphere~\cite{CuiFreeden}, and were obtained from the SpherePts software package~\cite{SpherePts}.  The number $n=924$ corresponds to a truncated div-free  vector spherical harmonic expansion of $\mu_0 = 42$, which is commonly used in scalar spherical harmonic tests where it is denoted as ``T42''~\cite{FornbergPiret2007}.  We only present results for the MQ kernel, but note that similar results were obtained for other smooth kernels.

\begin{figure}[htb]
\centering
\begin{tabular}{cc}
\includegraphics[width=0.47\textwidth]{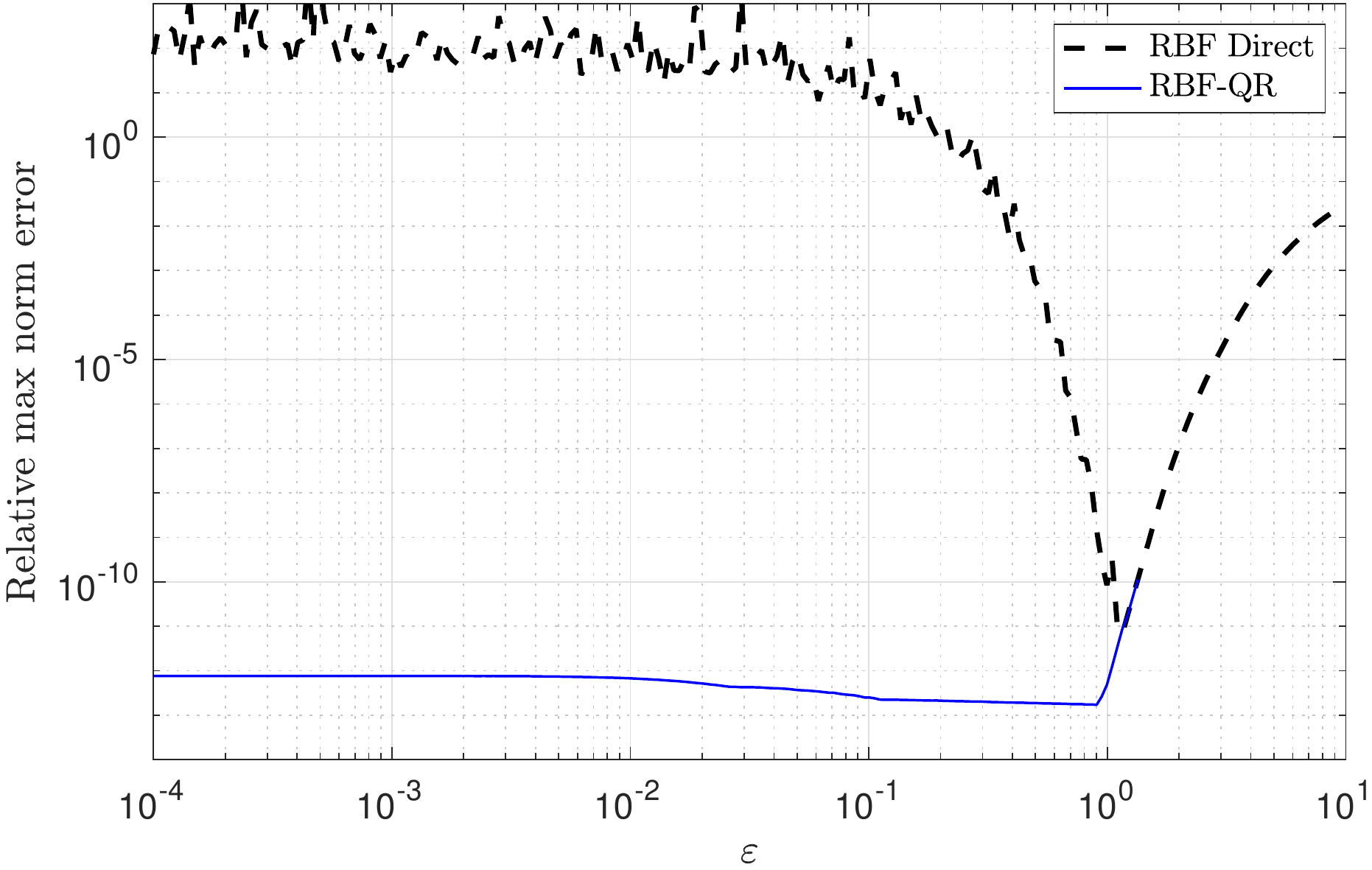}  &
\includegraphics[width=0.47\textwidth]{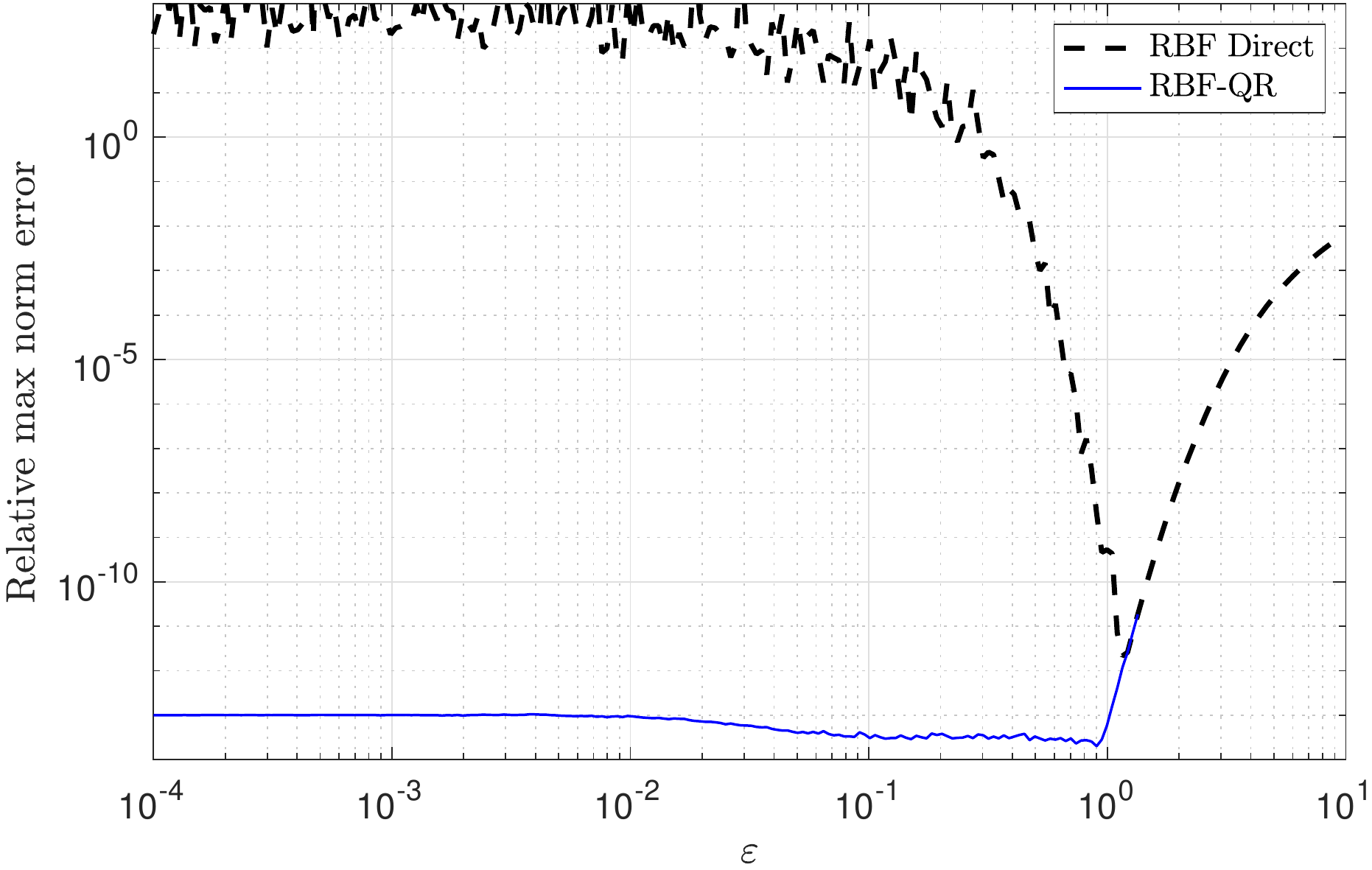} \\
(a) Field $\vu$ & (b) Stream function $\psi$
\end{tabular}
\caption{Comparison of the relative max-norm errors in reconstructing the target (a) div-free  field and (b) stream function \eqref{eq:target_sf} as a function of $\ep$ using the vector RBF-QR algorithm and RBF-Direct for the MQ kernel. \label{eq:errors_vs_ep}}
\end{figure}

Using this target field and node set, we computed the surface div-free  RBF interpolants using the RBF-Direct and vector RBF-QR algorithms for various values of $\ep$.  The relative max-norm error in the interpolant was then computed by evaluating the difference between the interpolants and the target field at a denser set of $4n$ points over the sphere.  The results from the experiment are displayed in Figure \ref{eq:errors_vs_ep} (a).  We see from the figure that the error in the RBF-Direct approach decreases rapidly with decreasing $\ep$ until around $\ep = 1$ when it starts to increase exponentially fast.  This is where ill-conditioning in the standard div-free  RBF basis sets in.  The vector RBF-QR algorithm on the other hand, shows no issues with ill-conditioning for any values of $\ep \leq 1$ and we can use it to compute the resulting interpolant in a stable manner all the way to the flat limit.  Note that the error reaches a minimum at a non-zero value of $\ep$ and then starts to increase slightly.  This is a feature that comes from the target field, and not from instabilities in the algorithm.  Figure \ref{eq:errors_vs_ep} (b) shows the max-norm errors in approximating the stream function \eqref{eq:target_sf} as a function of $\ep$.  For this computation, we used \eqref{eq:rbf_sf} to extract a stream function from the the interpolants and then adjusted it to have the same mean as the target stream function \eqref{eq:target_sf}.  We see that there is a similar trend in these results, which is to be expected since the stream function comes from the interpolant. 


%
%
\vspace{-0.1in}
\section{Concluding remarks}
The vector RBF-QR algorithm presented here can be used to bypass the ill-conditioning associated with surface div-free RBF interpolation on a sphere in the $\ep\rightarrow 0$ limit.  This allows for more comprehensive studies of how $\ep$ affects the accuracy of the interpolants.  The derivation of the algorithm additionally demonstrates the connection between these interpolants and div-free vector spherical harmonics in the flat limit.  The algorithm can also be applied straightforwardly to surface curl-free RBF interpolation, since in this case, one simply has to replace the matrix-valued kernel $\Psidiv$ with the kernel {$\Phi_{\rm curl} = (P({\vx})\nabla_\vx)(P({\vy})\nabla_\vy)^T\phi(\|\vx - \vy\|)$, where $P(\vxi) = (I-\vxi\vxi^T)$}~\cite{DrakeThesis}.  Div-free and curl-free RBFs are also available for interpolation in $\mathbb{R}^d$~\cite{fuselier08}.  For this problem it may be possible to extend the RBF-QR algorithms~\cite{FLF11,FaMC12} developed for scalar RBF interpolation in $\mathbb{R}^d$.  {However, the algorithm in~\cite{FLF11} is limited to the Gaussian kernel and the algorithm in~\cite{FaMC12} is limited to kernels with known Mercer series expansions.}

%
%
\vspace{-0.15in}
\section*{Acknowledgments}
KPD's work was partially supported by the SMART Scholarship, which is funded by The Under Secretary of Defense-Research and Engineering, National Defense Education Program / BA-1, Basic Research. GBW's work was partially supported by National Science Foundation grant CCF 1717556. 

\vspace{-0.15in}
\bibliographystyle{elsarticle-num}
\bibliography{refs}	

\end{document}